\documentclass[english, 11pt]{article}
\usepackage{amsthm}
\usepackage{amsmath}
\usepackage{amssymb}
\usepackage{cite}
\usepackage{latexsym}
\usepackage{verbatim}
\usepackage{wrapfig, lipsum}
\usepackage{color}

\usepackage{color}
\usepackage{graphicx,epsfig,subfigure,dcolumn,bm,mathrsfs}
\makeatletter
\theoremstyle{plain}

  \theoremstyle{plain}
  
  \theoremstyle{plain}
  
  \theoremstyle{plain}
  
 \theoremstyle{definition}
 \newtheorem*{defn*}{\protect\definitionname}
  \theoremstyle{remark}
  \newtheorem*{rem*}{\protect\remarkname}

\makeatother

\usepackage{babel}
  \providecommand{\corollaryname}{Corollary}
  \providecommand{\definitionname}{Definition}
  \providecommand{\lemmaname}{Lemma}
  \providecommand{\propositionname}{Proposition}
  \providecommand{\remarkname}{Remark}
\providecommand{\theoremname}{Theorem}

\begin{document}

\title{Noisy Hegselmann-Krause Systems: Phase Transition and the $2R$-Conjecture}
\author{
{\sc Chu Wang}
\thanks{The Program in Applied and Computational Mathematics, Princeton University, Princeton,
NJ 08540. ({\tt chuw}@{\tt math.princeton.edu}).
}
\and
  {\sc Qianxiao Li}
\thanks{The Program in Applied and Computational Mathematics, Princeton University, Princeton,
NJ 08540. ({\tt qianxiao}@{\tt math.princeton.edu}). This author is supported by the Agency for Science,
Technology and Research, Singapore. 
}
\and
{\sc Weinan E}
\thanks{Department of Mathematics and the Program in Applied and Computational Mathematics, Princeton University, Princeton,
NJ 08540. ({\tt   weinan}@{\tt math.princeton.edu})}
\and
{\sc Bernard Chazelle}
\thanks{Department of Computer Science, Princeton University, Princeton, NJ 08540.
({\tt   chazelle}@{\tt cs.princeton.edu}). This author was supported in part by NSF grants CCF-0963825,
CCF-1016250 and CCF-1420112.
}
}

\maketitle
\begin{abstract}
The classic Hegselmann-Krause ({\em HK}\,) model for opinion dynamics consists of a set of
agents on the real line, each one instructed to move, at every time step, 
to the mass center of all the agents within a fixed distance~$R$. 
In this work, we investigate the effects of noise in the continuous-time version of the model
as described by its mean-field limiting Fokker-Planck equation.
In the presence of a finite number of agents, 
the system exhibits a phase transition from order to disorder as the noise increases.
The ordered phase features clusters whose width depends only on the noise level.
We introduce an order parameter to track the phase transition and resolve the  
corresponding phase diagram.
The system undergoes a phase transition 
for small $R$ but none for larger $R$.
Based on the stability analysis of the mean-field equation,
we derive the existence of a forbidden zone for the disordered phase to emerge.
We also provide a theoretical explanation for the well-known 
$2R$ conjecture, which states that, for a random initial distribution in a fixed interval,
the final configuration consists of clusters separated by a distance of roughly $2R$.
Our theoretical analysis also confirms previous simulations 
and predicts properties of the noisy {\it HK} model in higher dimension.
\end{abstract}

\section{Introduction}

Network-based dynamical systems have received a surge of attention lately.
In these systems, typically, 
a set of agents will interact by communicating through a dynamic graph
that evolves endogenously. The popularity of the model derives from its
widespread use in the life and social sciences~\cite{axelrod1997complexity,blondel2009krause,
castellano2009statistical,chazelle2012dynamics,chazelle2015Algo,easley2012networks,jadbabaieLM03}.
Much of the difficulty in analyzing these systems
stems from the coupling between agent dynamics and evolving graph topology~\cite{chazelle2012dynamics}.
If the system is diffusive and the information transfer between agents 
is symmetric, it usually converges to an attractor 
under mild assumptions \cite{chazelle2011total, moreau2005stability}.
In the absence of symmetry, however, 
the system can exhibit the whole range of dynamical regimes,
from periodicity to chaos~\cite{chazelle2012dynamics}.

The Hegselmann-Krause ({\it HK}\,) model is the classic representative
of the diffusive type. It consists of a fixed number $N$ of agents, each one located at $x_k(t)$ on the real line.  
At each time step, every agent moves to the mass center of all the others within a fixed distance $R$. 
The position of an agent represents its ``opinion.''  The underlying assumption is that 
people are immune to the influence of others whose opinions greatly differ from their own.
In particular, two groups of agents initially separated by a distance of $R$ or more
will form decoupled dynamical systems with no interaction between them.
{\em HK} systems are known to converge in finite time, but
the relationship between the initial and final profiles remains mysterious. 
The celebrated {\em $2R$ conjecture} states, for a random initial distribution in a fixed interval,
the final configuration consists of clusters separated by a distance of roughly $2R$~\cite{blondel20072r}.

It is natural to enlarge the model by introducing noise into the dynamics~\cite{pineda2013noisy}. 
Stochasticity can be invoked to capture nonobservable factors operating at smaller scales.
Analytically, it also has the benefits of nudging the system away from pathological configurations.
By tuning the noise level as we would the temperature of a thermodynamical system,
we can vary the dynamics from chaos to fixed-point attraction and uncover
phase transitions in the process.  To simplify the analysis, we model
the system with a stochastic differential equation for the continuous-time
version of the {\em HK} model and focus on its mean-field approximation
in the form of a Fokker-Planck type partial differential equation governing the agent density evolution.
This formulation of noisy {\em HK} systems in the thermodynamic limit 
can be derived from first principles and seems well-supported by computer simulation.

\paragraph{Results and organization of the paper.}

After the formal introduction of the model in Section~\ref{model},
its long-time behavior is analyzed in Section~\ref{property}, along with a discussion
of basic path properties.  The bistability of the system is investigated 
in Section~\ref{phaseT} and 
an order parameter is introduced to track the phase transition and resolve the  
corresponding phase diagram.
We find that the system undergoes a phase transition 
for small $R$ but none for larger $R$.
Based on the stability analysis of the mean-field equation (Section~\ref{2R}),
we derive the existence of a forbidden zone for the disordered phase to emerge.
This puts us in a position to provide a theoretical explanation, 
to our knowledge the first of its kind, for the $2R$ conjecture.
Our theoretical analysis confirms previous simulations 
and predicts properties of the noisy {\it HK} model in higher dimension.
The pseudo-spectral scheme used for 
our simulations is discussed in the Appendix.

\paragraph{Prior work.}
 
The convergence of the classical {\it HK} system was established in a number of 
articles~\cite{hendrickx2006convergence,lorenz2005stabilization,moreau2005stability}
and subsequent work provided increasingly tighter bounds on the convergence rate,
with a current bound of $O(N^3)$~\cite{bhattacharyya2013convergence}.
While there exists a worst-case lower bound of $\Omega(N^2)$~\cite{DBLP:journals/corr/WedinH14a},
computer simulations suggest that, in practice, the convergence rate is at most linear.
The model extends naturally to higher dimension by interpreting $R$ as
the Euclidean distance, and polynomial bounds are known for that case
as well~\cite{chazelle2011total,bhattacharyya2013convergence}.
We note that the convergence time can be significantly lowered if
certain ``strategic" agents are allowed to move anywhere at each step~\cite{2014arXiv1411.4814K}.
For general consensus and stability properties of the infinite-horizon profile, we refer 
the interested readers
to~\cite{fortunato2005consensus,martinez2007synchronous,touri2011discrete,lorenz2006}.

Attempts to generalize the {\em HK} model to the nonsymmetric case have
proven surprisingly frustrating.  While it is known that diffusive influence systems
(the generalization of {\em HK} model) can have limit cycles and even chaotic
behaviors, the simple fact of allowing each agent to pick its own interval
produces dynamics that remains unresolved to this day. 
Numerical simulations suggest that such systems converge but a proof has been elusive.
All we know is that if each agent can pick its interval freely in $\{0,R\}$,
the system still converges\cite{DBLP:journals/corr/ChazelleW15};
in other words, taking the original {\it HK} system and fixing some of
the agents once does not change the fact that all the orbits have fixed-point attractors.

Regarding the $2R$ conjecture, 
the concept of equilibrium stability was introduced in \cite{blondel20072r} 
to put this conjecture on formal grounds.
Extensive experiments were conducted, suggesting a closer to $2.2R$.
All the work cited so far considers only the deterministic version of the system.
For the noisy version of the model,
Pineda {\it et al.} consider a discrete-time formulation where, at each step,
every agent randomly chooses to perform the usual {\it HK} step 
or to move randomly~\cite{pineda2013noisy}. 
Two types of random jumps are considered: bounded jumps confine agents
to a bounded distance from their current position while free jumps allow them to move anywhere.
An approximate density-based master equation is adopted for the analysis of the order-disorder phase transition and the noisy {\it HK} system is then compared with another opinion dynamics system,
the so-called {\em DW} model~\cite{deffuant2000mixing}.

\section{The Model \label{model}}

The stochastic differential equation (SDE) model we use in this paper can be expressed as
\begin{equation}\label{SDE}
\, dx_i=-\frac{1}{N}\sum_{j :\,|x_i-x_j|\le R}(x_i-x_j)\, dt+\sigma\, dW^{(i)}_t,
\end{equation}
where $i=1,2,\dots, N$ denotes the agents, 
$\sigma$ specifies the magnitude of the noise and
$W_t^{(i)}$ represent independent Wiener processes.
For technical convenience, we impose periodic 
boundary conditions on~(\ref{SDE}) by taking each $x_i$ modulo $L=1$
and interpreting $|z|$ as $\min \{|z|, L-|z| \}$.
Intuitively, the model mediates the tension between two competing forces:
the summation in~(\ref{SDE}) represents an attracting force that pulls 
the agents together while the diffusion term keeps them active in Brownian motion.
This can be compared to the use of two parameters in the noisy model of
Pineda {\it et al.}: a noise intensity $m$ determines the probability that
an agent should move to the mass center of its neighbors (vs. moving randomly),
and $\gamma$ bounds the length of the jump \cite{pineda2013noisy}. 
In the continuous-time model,
the pair $(m, \gamma)$ reduces to a single parameter, namely the noise magnitude $\sigma$.

\paragraph{The mean-field approximation.}

In the mean field limit $N\rightarrow \infty$, Equation~(\ref{SDE}) induces a nonlinear Fokker-Planck equation for the agent density profile $\rho(x,t)$~\cite{DBLP:journals/corr/GarnierPY15}:
\begin{equation}\label{PDE}
\rho_{t}(x,t)=\bigg{(}\rho(x,t)\int (x-y)\rho(y,t)\bm{1}_{|y-x|\le R})\, dy\bigg{)}_{x}+\frac{\sigma^2}{2} \rho_{xx}(x,t).
\end{equation}
The function $\rho(x,t)$ is the limiting density of
$\rho^N(x,t) := \frac{1}{N}\sum \delta_{x_j(t)}(dx)$,
as $N$ goes to infinity, where $\delta_x(dx)$ denotes
the Dirac measure with point mass at $x$.
In this partial differential equation form, the second derivative term represents the diffusion process that flattens the density $\rho$. On the other hand, the first term represents the advection of the density caused by attraction.
In higher dimensions, one just needs to replace the first derivative by a 
divergence operator and replace the second derivative by a Laplace operator.
We will use bold letter to denote a vector or a point in higher dimensional space.

To derive~(\ref{PDE}) from Equation (\ref{SDE}), we
consider the agent flux of a small piece of space $\Omega$. 
The Brownian motion of infinitely many agents is equivalent to a diffusion process,
which implies that the flux caused by the noise is $-\sigma\nabla \rho$.
On the other hand, the attraction between agents causes a flow with velocity 
$\int_{|{y}-{x}|\leq R}({y}-{x})\rho({y},t)\, d{y}$ at ${x}$;
hence a net outflow equal to the derivative of
$\rho({x},t)\int_{|{y}-{x}|\leq R}
      ({y}-{x})\rho({y},t)\, d{y}$.
Equation (\ref{PDE}) follows immediately from mass conservation and
the divergence theorem.
Regarding the boundary condition, none is necessary if we consider the system on the real line.
For the case of an interval, we use Neumann boundary conditions, in this case
a reflecting boundary condition since the flux at the boundary should be zero.
To simplify our analysis and simulation,  
we use a periodic boundary condition over the unit interval
for the rest of this paper. It should be noted that boundary effect plays a minor role in the dynamics of the system.

\section{Long-time Dynamics\label{property}}

As we mentioned earlier, the original {\it HK} system always 
converges within a  number of steps and the final configuration
consists of a union of clusters with pairwise distance larger than $R$.
A case of particular interest is that of a single cluster, forming what is
commonly called {\em consensus}. 
Of course, when noise is added to the system, there is no ``final'' state
to speak of; nevertheless, for the SDE model, one could focus 
on the averaged long-time behavior of the system as measured by the number of clusters.
Intuitively, a higher noise level $\sigma$ should correspond to more diffused clusters while,
for $\sigma$ above a certain threshold, clusters should break apart and
release the agents to move randomly in the unit interval.
This intuition is confirmed by simulating (\ref{SDE}) for different 
values of $R$ and $\sigma$.

\begin{figure}[htp!]
\includegraphics[width=\textwidth]{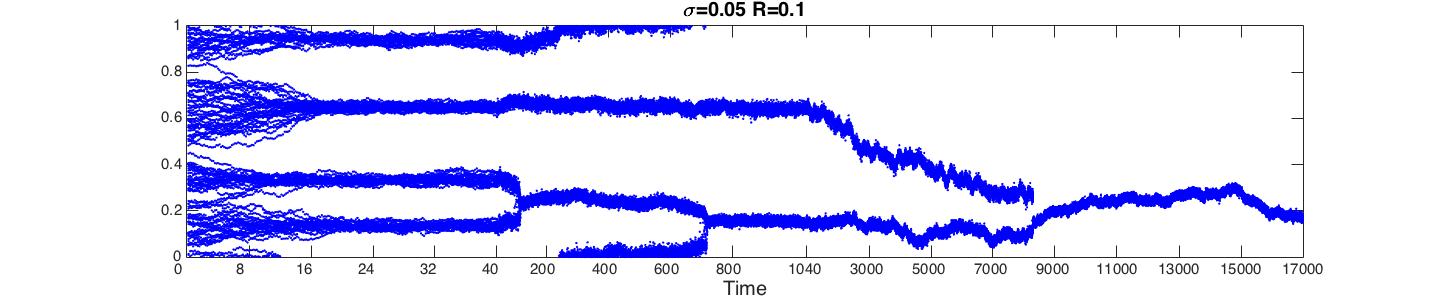}
\caption{Simulation of the SDE model for noisy {\it HK} system with $n=100$, $\sigma=0.05$ and $R=0.1$. Several clusters are formed at the beginning and later merge with each other.
The time axis is suitably transformed to make the evolution more apparent.
\label{path1}}
\end{figure}

\begin{figure}[htp!]
\includegraphics[width=\textwidth]{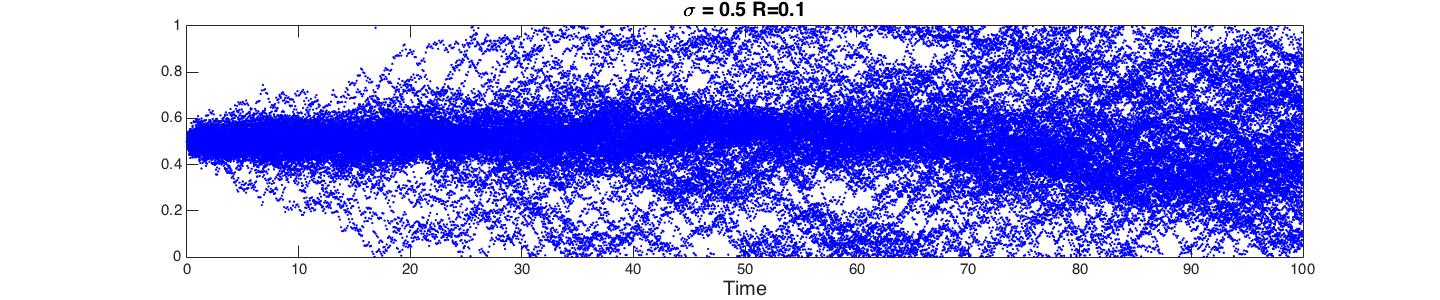}
\caption{Simulation of a noisy {\it HK} system with $n=100$, $\sigma=0.5$ and $R=0.1$.
All the agents start at $x=0.5$ and then diffuse to fill up the space.\label{path2}}
\end{figure}

We highlight two interesting scenarios in the SDE model.
For small values of~$\sigma$, a random initial distribution of the agents evolves 
in the following manner: at the beginning, the attracting forces dominate
and break symmetry by forming clusters of width less than $2R$.
As the clusters are formed, the noise term gradually overtakes the dynamics
and produces a jiggling motion of the clusters.
The mass center of the clusters follows 
a Brownian motion of variance $\sigma^2/n$, where $n$ is the number of agents in the cluster. 
Brownian motion in one (and two) dimension(s) is recurrent
so, as one would expect, the clusters formed in the early stage will eventually
merge almost surely. Such merging process can be observed in Figure~\ref{path1}
(with the time axis suitably transformed to make the evolution more apparent).
One does not expect this phenomenon in dimension three and above unless 
periodic conditions are added.
If  $\sigma$ is large and we initialize the system by placing all the agents at the same position, the 
initial cluster breaks apart to fill up the entire space (Figure~\ref{path2}).

\vspace{1cm}

\vspace{1cm}

\vspace{1cm}

\begin{figure}[htp!]
        \centering
        \begin{minipage}{0.32\textwidth}
                \centering
                \includegraphics[width=\textwidth]{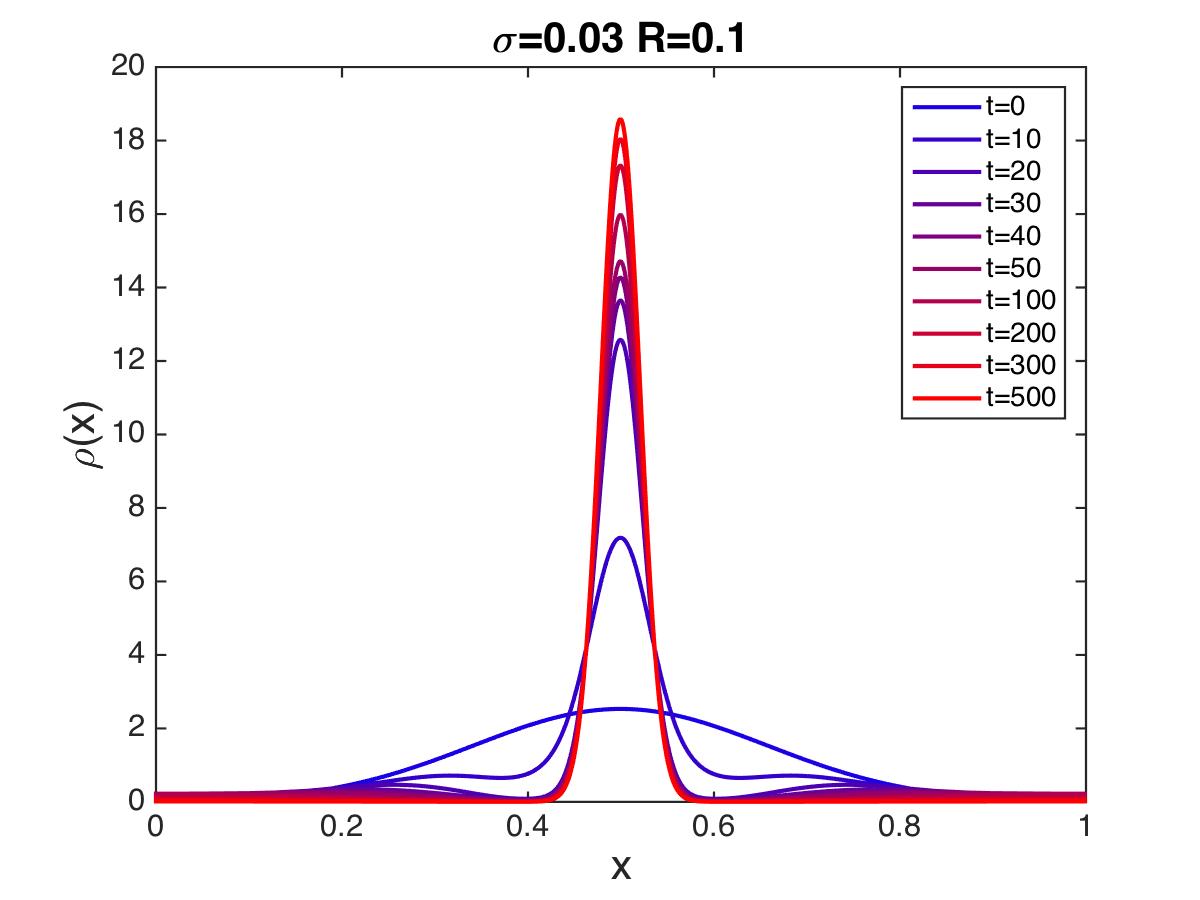}
                        \end{minipage}
      \begin{minipage}{0.32\textwidth}
                \centering
                \includegraphics[width=\textwidth]{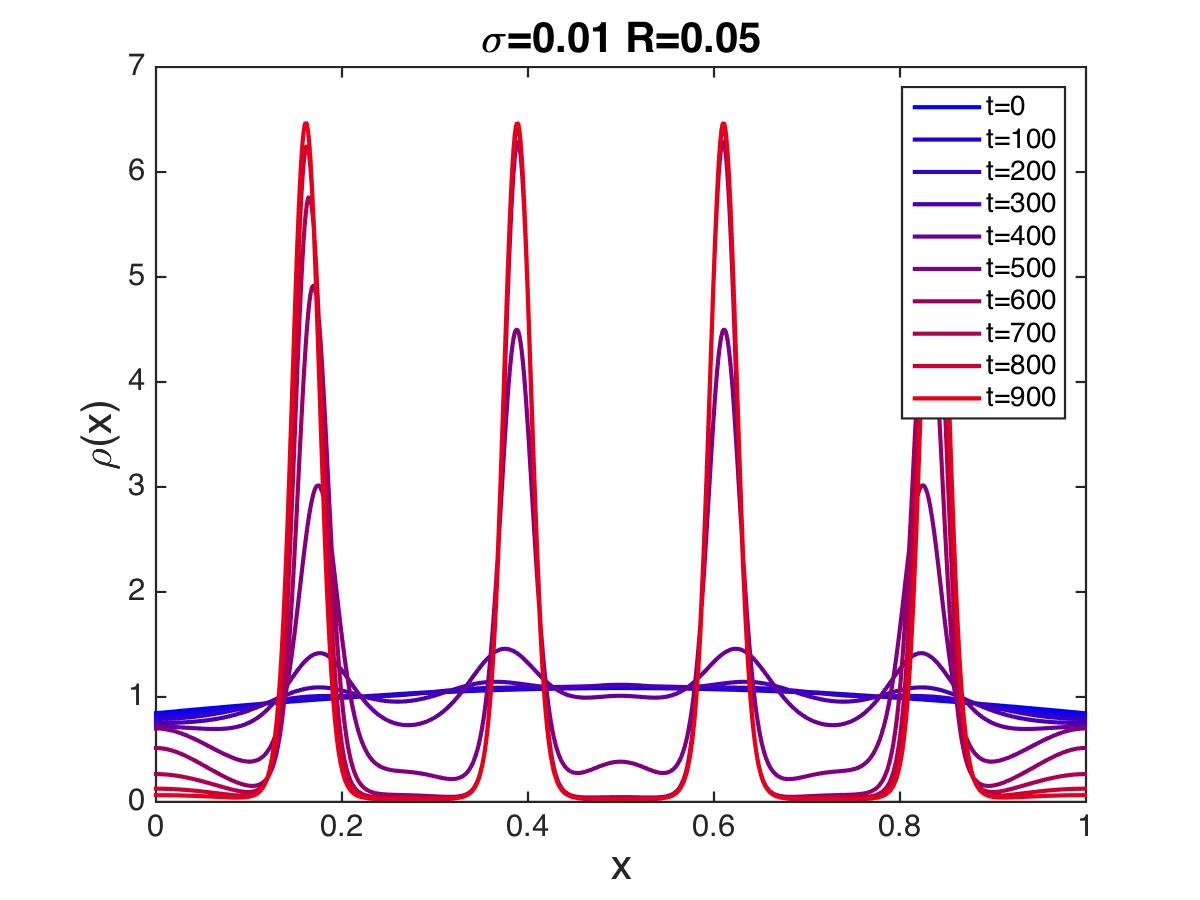}
        \end{minipage}
      \begin{minipage}{0.32\textwidth}
                \centering
                \includegraphics[width=\textwidth]{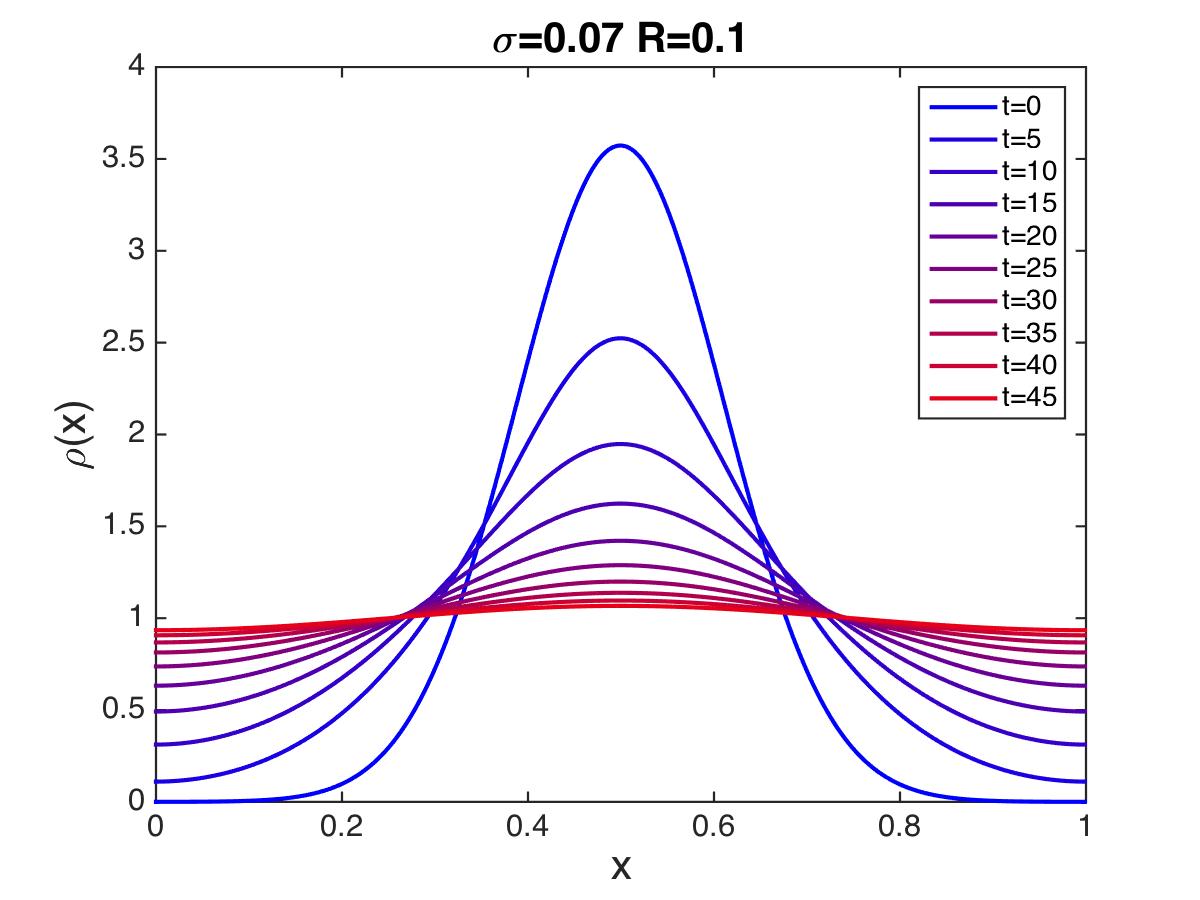}
        \end{minipage}
\caption{$\rho(x,0)\propto e^{-20(x-0.5)^2}$ (left);  $\rho(x,0)\sim e^{-(x-0.5)^2}$ (middle);
$\rho(x,0)\sim e^{-40(x-0.5)^2}$ (right). \label{simulation}}
\end{figure}

For the PDE model, we use a semi-implicit pseudo-spectral method~\cite{fornberg1998practical} for simulating (\ref{PDE}).
The details of the numerical procedure are given in the Appendix.  
In the simulation, the initial profile is rescaled before the system starts to 
ensure that the total mass is equal to 1. In Figure~\ref{simulation} (left and middle),
a relatively small $\sigma$ is chosen, with $R$ smaller in the middle figure,
which results in multiple clusters. 
In the right figure, even though the initial profile is a cluster, 
the density tends to flatten because of the high noise.

\vspace{1cm}

\section{Phase Transition and Order Parameter \label{phaseT}}

\begin{figure}[htp!]
        \centering
        \begin{minipage}{0.32\textwidth}
                \centering
                \includegraphics[width=\textwidth]{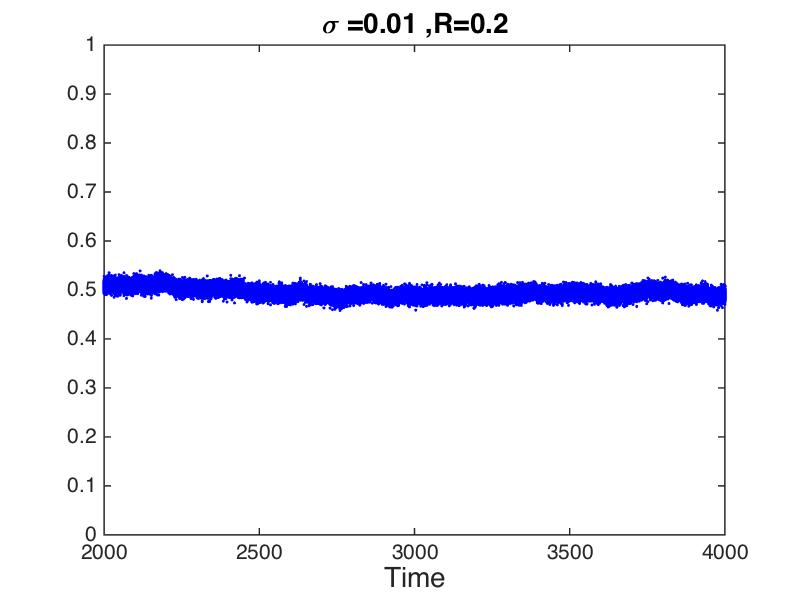}
                        \end{minipage}
      \begin{minipage}{0.32\textwidth}
                \centering
                \includegraphics[width=\textwidth]{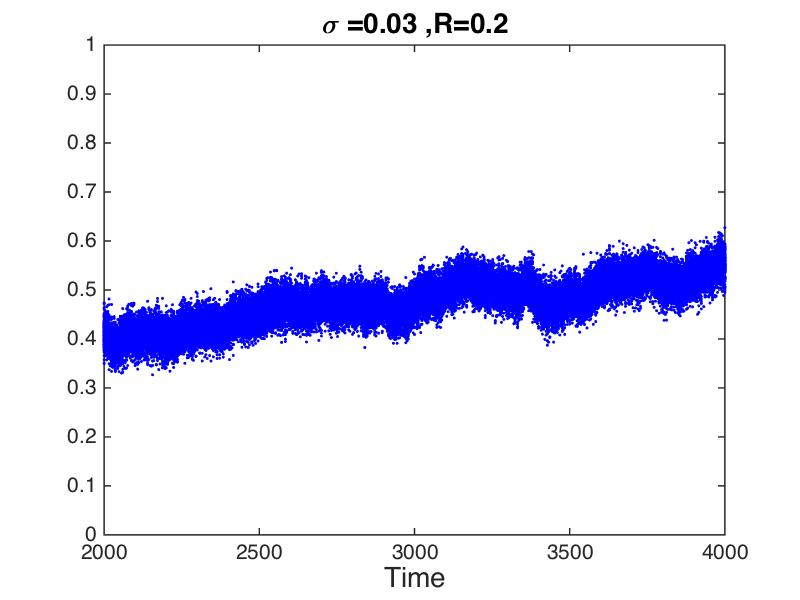}
        \end{minipage}
      \begin{minipage}{0.32\textwidth}
                \centering
                \includegraphics[width=\textwidth]{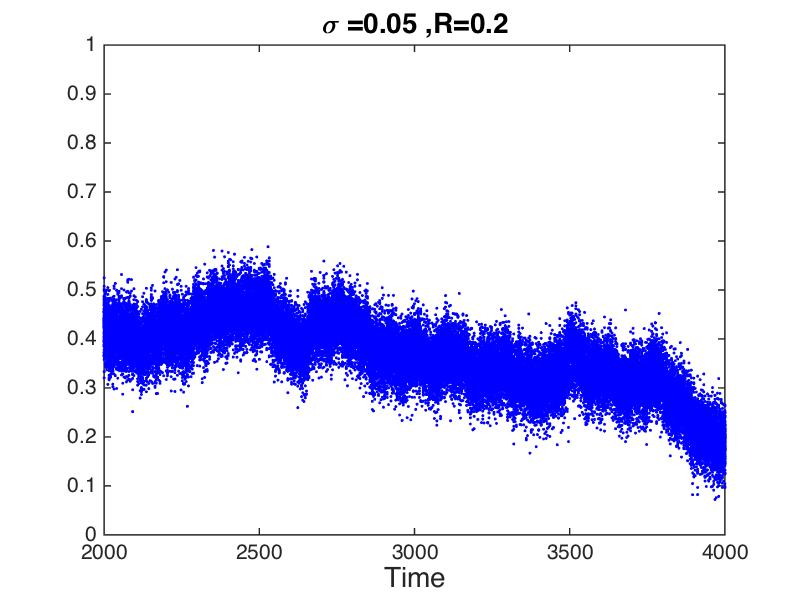}
        \end{minipage}
      \begin{minipage}{0.32\textwidth}
                \centering
                 \includegraphics[width=\textwidth]{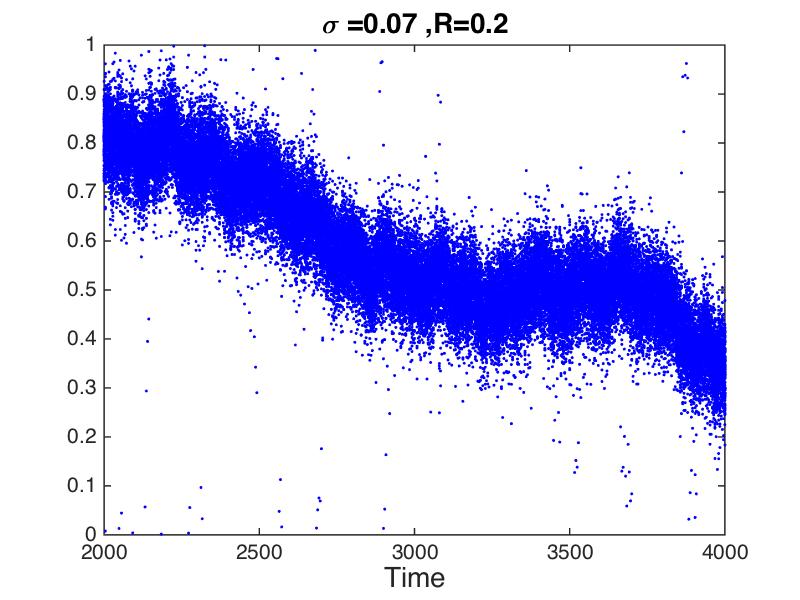}
        \end{minipage}      
        \begin{minipage}{0.32\textwidth}
                \centering
                \includegraphics[width=\textwidth]{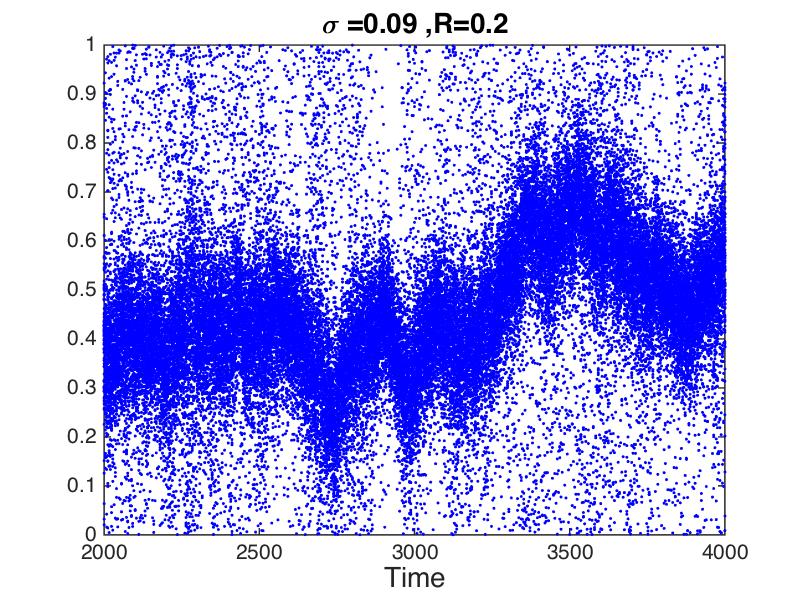}
        \end{minipage}
      \begin{minipage}{0.32\textwidth}
                \centering
                 \includegraphics[width=\textwidth]{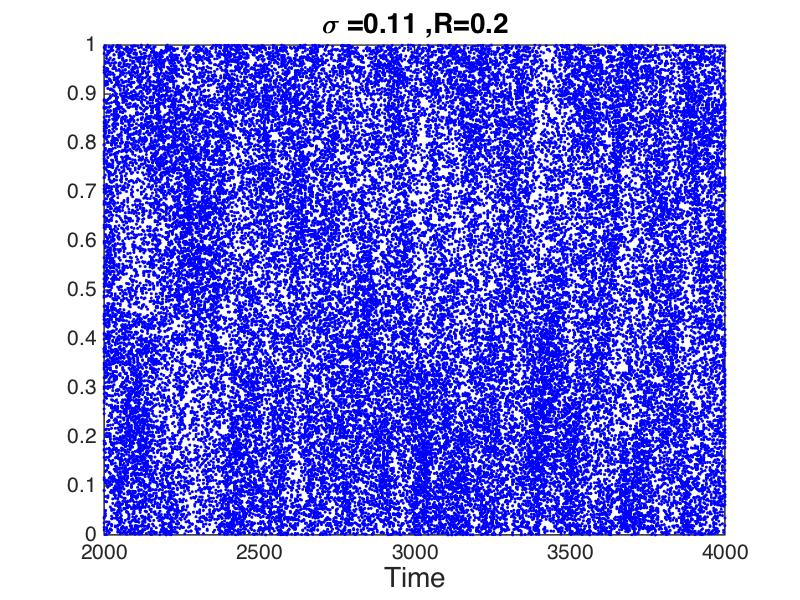}
        \end{minipage}
 \caption{Simulations of six independent systems,
 for $N=100$, $R=0.2$, and $\sigma=0.01, 0.03, \dots, 0.11$.
 The display is limited to the time window
 $[2000, 4000]$.  For visualizing purposes, the diagram is shifted
 to keep the cluster around the center.
 We note how the cluster widens as the noise level $\sigma$ grows,
 until the latter reaches a critical value signalling the start of the disordered phase. \label{6}}
 \end{figure}

A low noise level $\sigma$ 
corresponds to clustered phase in which a single cluster performs a random walk.
Conversely, a large value of $\sigma$ is associated
with a disordered phase, as the noise comes to dominate
the dynamics and render the attraction between agents negligible.
Our simulations reveal a bistable behavior.
As $\sigma$ increases and crosses a certain threshold, 
the cluster becomes dispersed and finally fills up space.
When $\sigma$ approaches the critical point, 
some agents may escape the main cluster and travel through space but eventually
return to the cluster. 
In general, clusters widen as $\sigma$ grows.
This observation will be explained in Section~\ref{clusterwidth}.
To illustrate the change of behavior as the noise level increases,
sample paths for fixed $R$ are shown in Figure \ref{6}
for $\sigma = 0.01, 0.03, \dots, 0.11$.
Since the long-time behavior of the system does not depend on the initial value,
we always use random initial values for the agents and wait for the system to stabilize.

To describe and distinguish between the clustered and disordered phases, 
we introduce the order parameter 
$$Q(\bm{x}):= \, \frac{1}{N^2}\sum_{i,j=1}^N \mathbf{1}_{|x_i-x_j|\le R}\, $$
to measure the edge density of the communication graph.  
Obviously, $Q=1$ when a consensus state is reached and it can be easily
checked that $Q=2R/L$ when the system is in a disordered state.

\begin{figure}[htp!]
\centering
\includegraphics[width=0.7\textwidth]{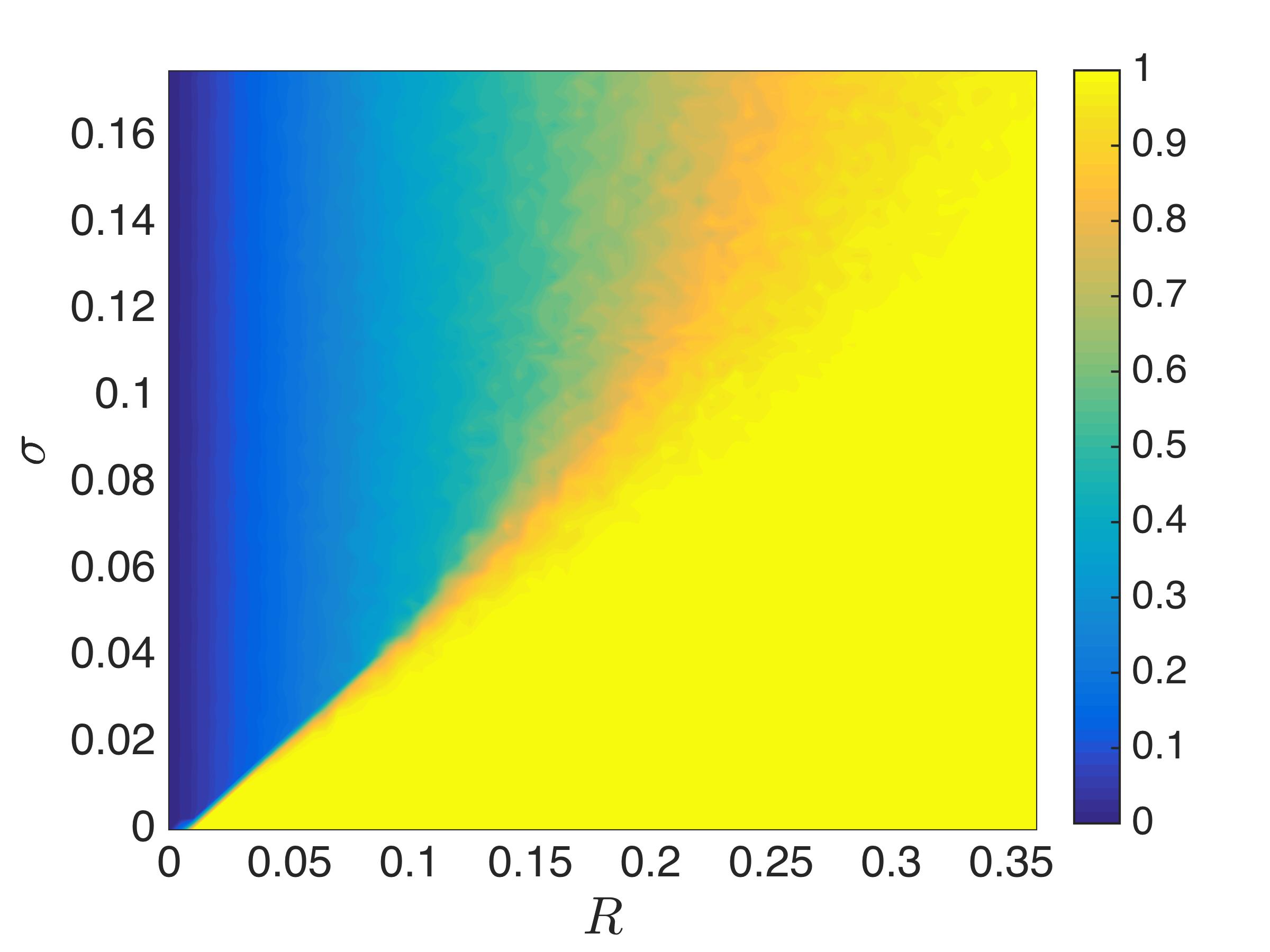}
\caption{Phase diagram\label{pd} for the SDE model. The color scale represents the value of $Q$
for $L=1$, $N=300$, and $t=10^5$, the latter being 
large enough to ensure that the system has reached steady state. 
We note the existence of a phase transition 
for small $\sigma$ consisting of a line separating clustered and disordered states. 
The transition line vanishes at higher noise level.}
\end{figure}

Figure \ref{pd} suggests that $Q(\bm{x})$ is discontinuous near the origin
but the discontinuity vanishes for large $R$.
In fact, when the interval length $R$ is large, 
the noise level $\sigma$ need to be comparably large to overcome
the attracting forces among the agents. 
The distinction between the clustered phase and the disordered
one becomes blurry for large $(R, \sigma$)
and the phase transition ceases to be observable.

\section{Analysis of the Clustered Phase\label{clusterwidth}}

We analyze the long-time profile of the clustered phase with respect to the interval length $R$
and the noise level $\sigma$.
An important difference with the 
original noiseless {\it HK} system is that clusters do not form single points
but evolving intervals of moving points.
We seek conditions on $(R,\sigma)$ to make the clusters stable.

\paragraph{The SDE model.}

If we assume that the widths of the clusters are smaller than $2R$, 
then the drift term becomes simpler.
Indeed, suppose there are $n$ agents in the cluster and 
let $g(t)=\frac{1}{n}\sum_{i=1}^nx_i(t)$ denote the center of mass.
The equations for $x_i(t)$ and $g(t)$ become
$\, dx_i=\frac{n}{N}(g(t)-x_i)\, dt+\sigma\, dW_t^{i}$ and
$\, dg=\frac{\sigma}{n}\sum_{i=1}^n\, dW_t^{(i)}$.
Then it is straightforward to see that 
\begin{equation}
\label{eq:sde_prof}
\, dx_i=-\frac{n}{N}x_i\, dt+\sigma\, dW_t^{i}+\frac{\sigma}{N}\sum_{k=1}^n\, dW_t^{k}.
\end{equation}
Pick one cluster and assume that all of its agents are initially placed at the origin.
They will oscillate in Brownian motion while being
pulled back to the origin because of the first term of~(\ref{eq:sde_prof}).
It follows that the invariant measure of the stochastic differential equation above 
provides a faithful description of the local profile. 
Since the SDE~\eqref{eq:sde_prof} is linear in $\bm{x}= (x_1,\ldots, x_n)$, it describes a Gaussian process;
therefore, it suffices to study the first two moments of the vector $\bm{x}$ 
to provide its complete characterization in distribution. 
Using It\^{o} calculus, we obtain the following moment equations:
\begin{eqnarray*}
\, d\mathbb{E}x_i & = & -\frac{n}{N}\, \mathbb{E}x_i \, dt \, ,\\
\, d\mathbb{E}x_i x_j & = & \left(-\frac{2n}{N} \, \mathbb{E}x_i x_j + \sigma^2(\Sigma\Sigma^T)_{ij}\right)\, dt,  
\end{eqnarray*}
where $\Sigma_{ij} = \frac{1}{N} + \delta_{ij}$. 
This implies that
$$\mathbb{E}x_i x_j  =  \frac{N\sigma^2}{2n} (\Sigma\Sigma^T)_{ij}(1-e^{-2nt/N}) .$$
Hence, at steady state, 
$$
\mathbb{E}x_i = 0  
\hspace{1cm} \text{and} \hspace{1cm} 
\mathbb{E}x_i x_j = \frac{N\sigma^2}{2n} (\Sigma\Sigma^T)_{ij}.  
$$
This implies a Gaussian profile at steady state of the form
$x_i\sim\mathcal{N}(0,\frac{N\sigma^2}{2n}\Sigma\Sigma^T)$. 
Notice that the covariance matrix $\frac{N\sigma^2}{2n}\Sigma\Sigma^T$ has one-fold eigenvalue $(1+\frac{n}{2N}+\frac{N}{2n})\sigma^2$
and $(n-1)$-fold eigenvalue $\frac{N\sigma^2}{2n}$.
Therefore the cluster size depends linearly on $\sigma$.
Furthermore, since these two eigenvalues both decrease when $n$ increases from 1 to $N$,
the cluster size actually increases as the number of agents declines.
This result may seem counterintuitive at a first glance because
more agents are naturally associated with larger clusters. 
In the SDE model, however,
more agents lead to stronger attraction to the mass center, which results in 
smaller-size clusters.

\paragraph{The PDE model.}

The clustered phase corresponds to a steady solution of~(\ref{PDE}).
Asymptotic analysis shows that the clusters are locally Gaussian.
Again, we assume that the clustered phase is stable for the given parameters $R$ and $\sigma$.
In addition, we focus on the case where $R$ is much smaller than $L=1$. 
By~(\ref{PDE}), the steady-state equation can be expressed as
\begin{equation}\label{eq:steady_state}
	\left(\rho(x)\int_{x-R}^{x+R}(y-x)\rho(y)\, dy\right)'=\frac{\sigma^2}{2}\rho''(x).
\end{equation}
After integrating the above equation once, we have
\begin{equation}\label{eq2}
\rho(x)\int_{x- R}^{x+R}(y-x)\rho(y)\, dy=\frac{\sigma^2}{2}\rho'(x)+C_1.
\end{equation}
Note that $\rho(x)\equiv 1$ is a solution for $C_1=0$. 
For a single-cluster profile, we may assume without loss of generality 
that the cluster is centered around~$0$ and the solution is symmetric:
$\rho(x)=\rho(-x)$.  By periodicity, we can confine our analysis to the interval
$[-\frac{1}{2},\frac{1}{2}]$, which implies $C_1=0$. 
Rearranging and integrating \eqref{eq2} once more yields
\begin{equation}\label{eq:rearranged}
\rho(x)=\rho(0)\exp \left\{ \frac{2}{\sigma^2}\int_{- 1/2}^{1/2}K(x,y)\rho\left(y\right)\, dy\right\},
\end{equation}
where 
\begin{equation}\label{eq:kernel_exp}
K\left(x,y\right)=\int_{0}^{x} \left(y-\xi\right) \mathbf{1}_{\left\vert y-\xi\right\vert\leq R}\, \, d\xi.
\end{equation}
It is easy to evaluate the kernel $K$.  For $\vert x\vert>2R$, we have
\begin{equation}\label{eq:kernel_soln}
K\left(x,y\right)= \begin{cases}
\frac{1}{2} (R+x-y) (R-x+y) & x-R\leq y\leq x+R,\\
\frac{1}{2} (y^2-R^2) & -R\leq y\leq R,\\
0 & \text{otherwise}. \\
\end{cases}
\end{equation}
For $0 \leq x\leq 2R$, we get
\begin{equation}\label{eq:kernel_soln2}
K\left(x,y\right)= \begin{cases}
\frac{1}{2} (R+x-y) (R-x+y) & R\leq y\leq x+R,\\
-\frac{1}{2} x (x-2y)  & x-R \leq y\leq R,\\
\frac{1}{2} (y^2-R^2) & -R\leq y\leq x-R,\\
0 & \text{otherwise}. \\
\end{cases}
\end{equation}
and for $-2R \leq x< 0$, we get
\begin{equation}\label{eq:kernel_soln3}
K\left(x,y\right)= \begin{cases}
\frac{1}{2} (R+x-y) (R-x+y) & x-R\leq y\leq -R,\\
-\frac{1}{2} x (x-2y)  & -R \leq y\leq x+R,\\
\frac{1}{2} (y^2-R^2) & x+R\leq y\leq R,\\
0 & \text{otherwise}. \\
\end{cases}
\end{equation}
Now, it is not easy to solve the integral functional equation \eqref{eq:rearranged} for $\rho$. Instead, we try to find an asymptotic solution for small $\sigma$. Concretely, we define
\begin{equation}\label{eq:rho_0}
\rho_0\left(x\right)=
C e^{ {- \min\{x^2,R^2\}/\sigma^2}  }.
\end{equation}
where the normalization constant $C$ ensures summation to 1.
We now show that $\rho_0$ solves \eqref{eq:rearranged}
up to the leading term in the expansion in $\sigma$.
Observe that, for $\sigma$ much smaller than $R$, the profile $\rho_0$ is concentrated around $x=0$, 
so we only need
to evaluate the integral in \eqref{eq:rearranged} near $y=0$ and thus ignore
error terms exponentially small in $\sigma$. 
For $\vert x\vert \leq R$, by~(\ref{eq:kernel_soln2}, \ref{eq:kernel_soln3}), we see that near $y=0$, 
\begin{equation}
K\left(x,y\right)=-\frac{1}{2} x (x-2y). 
\end{equation}
The exponent on the right hand side of \eqref{eq:rearranged} becomes
\begin{eqnarray*}
\frac{2}{\sigma^2}\int_{-1/2}^{1/2}K\left(x,y\right)\rho_0\left(y\right)\, dy 
&=& -\frac{1}{\sigma^2}\int_{-1/2}^{1/2}  x (x-2y)\rho_0\left(y\right)\, dy \\
&=& -\frac{x^2}{\sigma^2}\int_{-1/2}^{1/2}\rho_0(y) \, dy
= - \frac{x^2}{\sigma^2},
\end{eqnarray*}
which satisfies~\eqref{eq:rearranged}.
On the other hand, for $\vert x\vert > R$, near $y=0$, we have
\begin{equation}
K\left(x,y\right)=\frac{1}{2} (y^2-R^2).
\end{equation}
Thus, 
$$
\frac{2}{\sigma^2}\int_{-1/2}^{1/2}K\left(x,y\right)\rho_0\left(y\right)\, dy =
\frac{1}{\sigma^2}\int_{-1/2}^{1/2} (y^2-R^2)\rho_0\left(y\right)\, dy 
= - \frac{R^2}{\sigma^2} + O\left(\sigma\right). 
$$
Once again, \eqref{eq:rearranged} is satisfied, this time up to the leading term in
the expansion in $\sigma$.
We have thus established that, in the presence of a low level of noise,
the single-cluster steady state has a Gaussian profile with variance $\sigma^2/2$ 
near the cluster with exponential error decay as function of $\sigma$.
Notice that convex combinations of cluster profiles 
of the form~(\ref{eq:rho_0}) but centered at different locations in  (\ref{eq:steady_state})
generate error terms in $O(\sigma)$ as long as the distances between clusters are large enough to eliminate 
inter-cluster interaction.  As a result, multiple Gaussians are also steady-state solutions 
in the $\sigma$-error sense as long as the Gaussians are well-separated.
The asymptotic solution \eqref{eq:rho_0} matches remarkably well the simulation results 
for the PDE model given in Section~\ref{property}.

\section{The Disordered Phase and the $2R$-Conjecture\label{2R}}

\begin{figure}[htp!]
\centering
\includegraphics[width=0.8\textwidth]{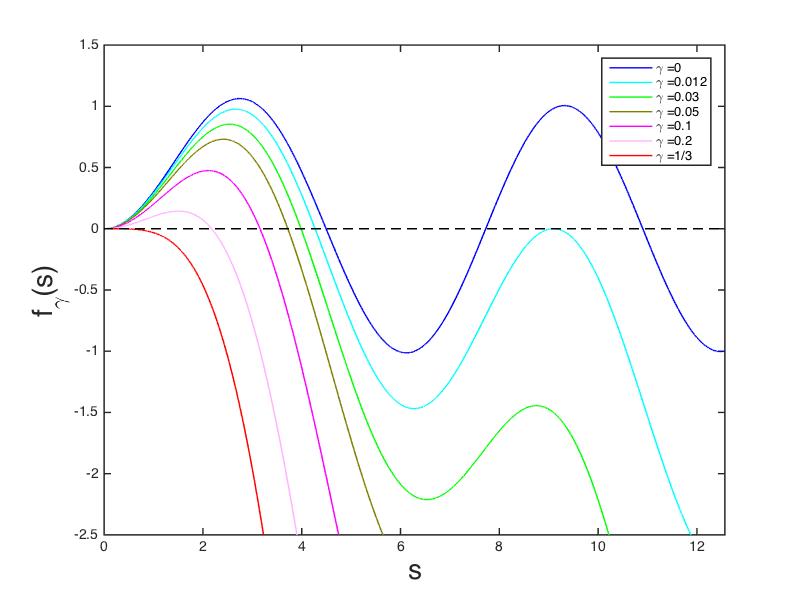}
\caption{The graph of $f_\gamma(s)$  for different values of $\gamma$.\label{f(s)}}
\end{figure}

In order to analyze the stability of the disordered phase,
we perturb the constant solution $\rho=1$ by Fourier waves and focus on the most unstable mode.
To unify the calculation in each dimension, we consider the mean-field
Fokker-Planck equation in dimension $d$ with periodic boundary conditions in the unit cube
(that is, reducing each coordinate modulo $L=1$):

\begin{equation}
\frac{\partial }{\partial t}\rho(\bm{x},t)=\nabla\cdot\left (\rho(\bm{x},t)
\int_{\|\bm{x}- \bm{y}\|\leq R}(\bm{x}-\bm{y})\rho(\bm{y},t)\, d\bm{y}\right )+\frac{\sigma^2}{2} \Delta\rho(\bm{x},t),
\end{equation}
$$
\rho(\bm{x},0)=\rho_0(\bm{x}),~~~~\int_{0}^{1}\rho(\bm{x},t)\, d\bm{x}=1,
~~~~0\le x_i\le 1,~~~~i=1,2,\dots,d.
$$
We perturb the constant solution by $\rho(\bm{x},t)=1+p(t)e^{2\pi i\bm{k}\cdot\bm{x}}$,
where $\bm{k}=(k_1,k_2,\dots,k_d)$. 
Gram-Schmidt orthogonalization guarantees that we could find an orthogonal matrix $M$ such that
the first row of $M$ is $\bm{k}/K$, where $K:=\|\bm{k}\|_2$.
By defining $\bm{z}=M\bm{y}$, we then have $z_1=\bm{k}\cdot\bm{y}/K$. 
Now let $s=2\pi K R$ and plug in the ansatz $\rho(\bm{x},t)=1+p(t)e^{2\pi i\bm{k}\cdot\bm{x}}$.
For $p$ small enough, this gives us
\begin{eqnarray*}
&p_te^{2\pi i\bm{k}\cdot\bm{x}}&\\
&=&
- \nabla\cdot\bigg{(}(1+pe^{2\pi i\bm{k}\cdot\bm{x}})
\int_{\|\bm{y}\|\le R}\bm{y}pe^{2\pi i\bm{k}\cdot(\bm{x}+\bm{y})}\, d\bm{y}\bigg{)} \\
&&  \hspace{5.5cm}  - 2\pi^2\sigma^2K^2pe^{2\pi i\bm{k}\cdot\bm{x}}\\
&\approx&
-2\pi i pe^{2\pi i\bm{k}\cdot\bm{x}}\int_{\|\bm{y}\|\le R}\bm{k}\cdot \bm{y}
e^{2\pi i\bm{k}\cdot\bm{y}}\, d\bm{y}
- 2\pi^2\sigma^2K^2pe^{2\pi i\bm{k}\cdot\bm{x}}\\
&\approx& -2\pi i pe^{2\pi i\bm{k}\cdot\bm{x}}\int_{\|\bm{z}\|\le R}Kz_1
e^{2\pi iKz_1}\, d\bm{z}
- 2\pi^2\sigma^2K^2pe^{2\pi i\bm{k}\cdot\bm{x}}\\
&\approx& s Rpe^{2\pi i\bm{k}\cdot\bm{x}}\int_{\|\bm{z}\|\le 1}z_1
e^{i s z_1} \, d\bm{z}
- 2\pi^2\sigma^2K^2pe^{2\pi i\bm{k}\cdot\bm{x}}\ .
\end{eqnarray*}
This results in the ODE
\begin{eqnarray*}
p_t&= & R\left (s\int_{\|\bm{z}\|\le 1}z_1\sin(sz_1)\, d\bm{z}
- \frac{\sigma^2}{2R^3}s^2\right )p\\
&=& 2Rp\left ( \frac{\sin s}{s} - \cos s - \frac{\sigma^2 }{4R^3} s^2\right ).
\end{eqnarray*}
In dimension $d=1$, 
the ODE reduces to $p_t/p=2R f_\gamma(s)$, 
where $\gamma:=\sigma^2/4R^3$ and $f_\gamma(s)$ is defined as
\begin{equation}\label{fs}
f_\gamma(s)=   \frac{\sin s}{s}  -\cos s - \gamma s^2.
\end{equation}
When $f_\gamma(s)\le 0$ for all $s=2k\pi R$, small perturbation to $\rho=1$ will decay and finally vanish.
On the other hand, if some $k$ makes $f_\gamma(s)>0$, then $\rho=1$ is no longer stable 
and the system yields a clustered phase. 
The function $f_\gamma(s)$ is graphed in Figure~ \ref{f(s)} for several values of $\gamma$.
Since its first two terms in~(\ref{fs}) are bounded, 
$f_\gamma(s)<0$ for all $s>s_0$ and $\gamma>0$, which means
that the high frequency modes are all stable. 
On the other hand, when $\gamma=0$ (the noiseless model),
$f_\gamma(s)$ behaves like $-\cos s$ for $s$ large enough,
which implies that $f_\gamma(s)$ can be positive for infinitely many frequencies.
When $\gamma>0.012$, $f_\gamma(s)>0$ only over an interval of the form $[0,s_1]$,
corresponding to the low frequency modess.
As $\gamma$ increases, $s_1$ shrinks to become nearly 0 at $\gamma=\frac{1}{3}$.

\vspace{1cm}

\paragraph{The unstable zone for the disordered phase.}

When the noise level $\sigma$ is very small,
the clustering effect of the system dominates and the system falls within the clustered phase.
A Taylor expansion of $f_\gamma(s)$ shows that 
it can take on positive values as long as $\gamma<\frac{1}{3}$, 
which gives us the critical curve $\sigma^2=\frac{4}{3}R^3$.
Below this curve lies the unstable zone for the disordered phase in which
symmetry breaking fragments the constant solution $\rho=1$ 
into separate clusters.

It should be noted that this boundary is accurate when $R$ is small enough,
since a suitable frequency $k$ can always be found 
such that $f_\gamma(s)>0$, hence making the disordered phase unstable.
Conversely, for large $R$,  $f_\gamma(s)$ 
could be negative for all $k$, conferring stability,
even though $f_\gamma(s)>0$  for some $s<2\pi R$.

\begin{figure}[htp!]
\centering
\includegraphics[width=0.7\textwidth]{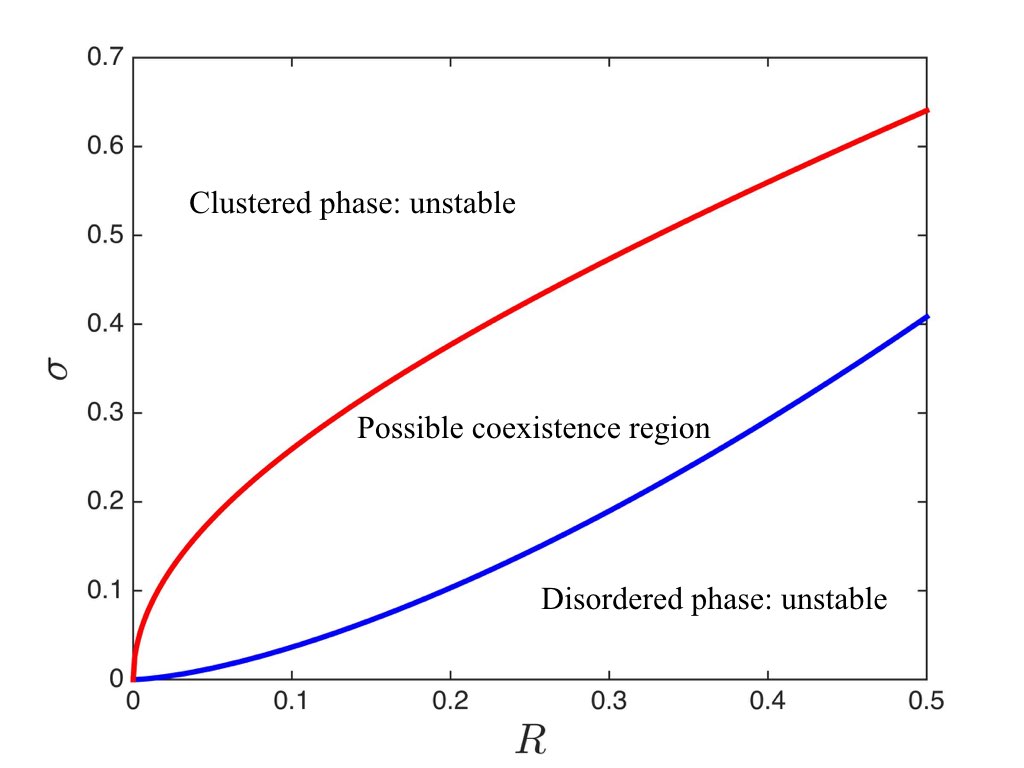}
\caption{Unstable zones for the disordered and clustered phases\label{zones}}
\end{figure}

Global stability conditions were discovered in \cite{2015arXiv151006487C},
showing that an unstable region for clustered phase lies
above the curve $\sigma^2=2(R+R^2/\sqrt{3})/\pi$.
Combining these two results, we see that
the two phases can coexist between the two curves, which
agrees with our simulations of the SDE---see (\ref{SDE})
and Figures~\ref{pd} and~\ref{zones}.

\vspace{1cm}

\paragraph{Resolving the $2R$ conjecture.}

The conjecture in questions refers to the original (noiseless) {\it HK} model:
it states that, if the agents starts out uniformly distributed in $[0, 1]$, 
they will converge toward clusters separated by distance of
roughly $2R$, thus setting their number close to the value $1/2R$.
This assumes that $N$ and $1/R$ are large enough.
To address this conjecture in the PDE model, we 
must confine $\gamma$ to the range $[0, \frac{1}{3}]$ to
make the clustered phase stable and then ask how many clusters one should expect.

We first rule out the case $s>2\pi$ by considering the pairwise distance of the resulting clusters.
Notice that when $s>2\pi$, $k$ is then larger than $1/R$. 
Even if the wave $e^{2\pi ikx}$ grows, it cannot exist for too long since the pairwise distance of the clusters would be smaller than $R$.
Therefore for this class of waves, no matter whether they are initially stable or unstable during the perturbation, they will not contribute to the final profile of the system. 
Then it is safe to focus on the interval $s\in [0,2\pi]$.
The expected number of clusters is given by the value of $s$ that maximizes
$f_\gamma(s)$.
As shown in Figure~\ref{f(s)}, the function increases to some peak and
then dips, while the root of $f_\gamma'(s)=0$ (the peak's location)
shifts to the left as $\gamma$ grows from $0$ to~$\frac{1}{3}$.
This implies that the number of clusters $k= s/(2\pi R)$ decreases as 
the noise level rises.
The $2R$ conjecture corresponds to the case $\gamma=0$.
Numerical calculation indicates that
the smallest nonzero root of $f_0'(s)$ is $s^*=2.7437$.
It follows that the expected number of final clusters is equal to
$$k^*=\frac{s^*}{2\pi R} =  \frac{1}{2.29R} \, .$$
Note that the bound comes fairly close to the number
observed experimentally for {\em HK} systems with a finite number of agents~\cite{blondel20072r}.

\paragraph{The higher-dimensional case.}

In higher dimensions, we have to consider functions of the form
$$F_\gamma(s):=\frac{1}{2}s\int_{\|\bm{z}\|\le 1}z_1\sin(sz_1)\, d\bm{z}
-\gamma s^2.$$
For the forbidden zone of the disordered phase, 
we can use a Taylor expansion for $F_\gamma$ in $s$ to determine the critical 
noise level $\sigma$ for a given $R$.
Let $S_{d-1}$ denote the area of the unit sphere in dimension $d$.
It is well known that $S_1=2\pi$ and $S_2=4\pi$.
For general $d$, we have 
$$S_{d-1}=\frac{2\pi^{d/2}}{\Gamma(\frac{d}{2})},~~\mathrm{where}~~\Gamma(s):=
\int_0^{\infty}x^{s-1}e^{-x}\, dx.$$
It follows that
\begin{eqnarray*}
F_\gamma(s) &\approx& \frac{1}{2}s\int_{\|\bm{z}\|\le 1}sz_1^2\, d\bm{z}-\gamma s^2\\
&\approx&\frac{s^2}{2d}\int_{\|\bm{z}\|\le 1} \|\bm{z}\|^2\, d\bm{z}-\gamma s^2\\
&\approx&\frac{s^2}{2d}\int_0^1 S_{d-1}r^{d+1} \, dr-\gamma s^2\\
&\approx&\left ( \frac{\pi^{d/2}}{d(d+2)\Gamma(\frac{d}{2})} -\gamma\right )s^2.
\end{eqnarray*}
The above expansion gives us the boundary at which the disordered phase loses stability:
$$\sigma^2=\frac{4\pi^{d/2}}{d(d+2)\Gamma(\frac{d}{2})}R^3.$$
Notice that the right-hand side equals $\frac{\pi}{2}R^3 > \frac{4}{3}R^3$ when $d=2$,
which means that in two dimensions we need even larger noise to make the disordered phase stable.
As the dimension grows, however, the Gamma function will bring the noise level down to zero
and the disordered phase will dominate unless the system becomes essentially noiseless.

\section{Conclusions\label{con}}

The contribution of this paper is
an analysis of the clustering modes of the 
noisy Hegselmann-Krause model.
We have provided theoretical insights, validated by numerical
simulations, into the stochastic differential equation model for 
a finite number of agents and the Fokker-Planck model for the 
mean-field approximation in the thermodynamic limit.
In the SDE model, we have shown that
the system exhibits either disorder or single-cluster profile.
In higher dimensions we conjecture that multiple clusters might be able to
co-exist. We have proposed an order parameter based on the edges density
of the communication graph to describe the phase transition.
In the PDE model, we used linear stability analysis to find a forbidden zone for disordered phase in which the constant solution cannot survive in the long run.
Most important, we provided a theoretical explanation,
the first to our knowledge, for the $2R$ conjecture.


\footnotesize{
\bibliography{hk}

\begin{thebibliography}{10}

\bibitem{axelrod1997complexity}
Robert Axelrod.
\newblock {\em The evolution of cooperation}.
\newblock Basic Books, 2006.

\bibitem{blondel2009krause}
Vincent~D. Blondel, Julien~M. Hendrickx, and John~N. Tsitsiklis.
\newblock On {K}rause's multi-agent consensus model with state-dependent
  connectivity.
\newblock {\em Automatic Control, IEEE Transactions on}, 54(11):2586--2597,
  2009.

\bibitem{castellano2009statistical}
Claudio Castellano, Santo Fortunato, and Vittorio Loreto.
\newblock Statistical physics of social dynamics.
\newblock {\em Reviews of modern physics}, 81(2):591, 2009.

\bibitem{chazelle2012dynamics}
Bernard Chazelle.
\newblock The dynamics of influence systems.
\newblock In {\em Foundations of Computer Science (FOCS), 2012 IEEE 53rd Annual
  Symposium on}, pages 311--320. IEEE, 2012.

\bibitem{chazelle2015Algo}
Bernard Chazelle.
\newblock An algorithmic approach to collective behavior.
\newblock {\em Journal of Statistical Physics}, 158:514--548, 2015.

\bibitem{easley2012networks}
David Easley and Jon Kleinberg.
\newblock {\em Networks, Crowds, and Markets: Reasoning About a Highly
  Connected World}.
\newblock Cambridge Univ. Press, 2010.

\bibitem{jadbabaieLM03}
Ali Jadbabaie, Jie Lin, and A.~Stephen Morse.
\newblock Coordination of groups of mobile autonomous agents using nearest
  neighbor rules.
\newblock {\em IEEE Trans. Automatic Control}, 48:988--1001, 2003.

\bibitem{chazelle2011total}
Bernard Chazelle.
\newblock The total s-energy of a multiagent system.
\newblock {\em SIAM Journal on Control and Optimization}, 49(4):1680--1706,
  2011.

\bibitem{moreau2005stability}
Luc Moreau.
\newblock Stability of multiagent systems with time-dependent communication
  links.
\newblock {\em Automatic Control, IEEE Transactions on}, 50(2):169--182, 2005.

\bibitem{blondel20072r}
Vincent~D Blondel, Julien~M Hendrickx, and John~N Tsitsiklis.
\newblock On the 2r conjecture for multi-agent systems.
\newblock In {\em Control Conference (ECC), 2007 European}, pages 874--881.
  IEEE, 2007.

\bibitem{pineda2013noisy}
Miguel Pineda, Ra{\'u}l Toral, and Emilio Hern{\'a}ndez-Garc{\'\i}a.
\newblock The noisy {H}egselmann-{K}rause model for opinion dynamics.
\newblock {\em The European Physical Journal B}, 86(12):1--10, 2013.

\bibitem{hendrickx2006convergence}
Julien Hendrickx and Vincent Blondel.
\newblock Convergence of different linear and non-linear {V}icsek models.
\newblock In {\em Proc. 17th Int. Symp. Math. Theory Networks Syst. (MTNS
  2006)}, pages 1229--1240, 2006.

\bibitem{lorenz2005stabilization}
Jan Lorenz.
\newblock A stabilization theorem for dynamics of continuous opinions.
\newblock {\em Physica A: Statistical Mechanics and its Applications},
  355(1):217--223, 2005.

\bibitem{bhattacharyya2013convergence}
Arnab Bhattacharyya, Mark Braverman, Bernard Chazelle, and Huy~L Nguyen.
\newblock On the convergence of the hegselmann-krause system.
\newblock In {\em Proceedings of the 4th conference on Innovations in
  Theoretical Computer Science}, pages 61--66. ACM, 2013.

\bibitem{DBLP:journals/corr/WedinH14a}
Edvin Wedin and Peter Hegarty.
\newblock A quadratic lower bound for the convergence rate in the
  one-dimensional hegselmann-krause bounded confidence dynamics.
\newblock {\em CoRR}, abs/1406.0769, 2014.

\bibitem{2014arXiv1411.4814K}
S.~{Kurz}.
\newblock {Optimal control of the convergence time in the Hegselmann--Krause
  dynamics}.
\newblock {\em ArXiv e-prints}, November 2014.

\bibitem{fortunato2005consensus}
Santo Fortunato.
\newblock On the consensus threshold for the opinion dynamics of
  {K}rause-{H}egselmann.
\newblock {\em International Journal of Modern Physics C}, 16(02):259--270,
  2005.

\bibitem{martinez2007synchronous}
Sonia Mart{\'\i}nez, Francesco Bullo, Jorge Cort{\'e}s, and Emilio Frazzoli.
\newblock On synchronous robotic networks -- part ii: Time complexity of
  rendezvous and deployment algorithms.
\newblock {\em Automatic Control, IEEE Transactions on}, 52(12):2214--2226,
  2007.

\bibitem{touri2011discrete}
Behrouz Touri and Angelia Nedi\'{c}.
\newblock Discrete-time opinion dynamics.
\newblock In {\em 2011 Conference Record of the Forty-Fifth Asilomar Conference
  on Signals, Systems and Computers (ASILOMAR)}, 2011.

\bibitem{lorenz2006}
Jan Lorenz.
\newblock Consensus strikes back in the hegselmann-krause model of continuous
  opinion dynamics under bounded confidence.
\newblock {\em Journal of Artificial Societies and Social Simulation}, 9(1):8,
  2006.

\bibitem{DBLP:journals/corr/ChazelleW15}
Bernard Chazelle and Chu Wang.
\newblock Inertial {H}egselmann-{K}rause systems.
\newblock {\em CoRR}, abs/1502.03332, 2015.

\bibitem{deffuant2000mixing}
Guillaume Deffuant, David Neau, Frederic Amblard, and G{\'e}rard Weisbuch.
\newblock Mixing beliefs among interacting agents.
\newblock {\em Advances in Complex Systems}, 3(01n04):87--98, 2000.

\bibitem{DBLP:journals/corr/GarnierPY15}
Josselin Garnier, George Papanicolaou, and Tzu{-}Wei Yang.
\newblock Consensus convergence with stochastic effects.
\newblock {\em CoRR}, abs/1508.07313, 2015.

\bibitem{fornberg1998practical}
Bengt Fornberg.
\newblock {\em A practical guide to pseudospectral methods}, volume~1.
\newblock Cambridge university press, 1998.

\bibitem{2015arXiv151006487C}
B.~{Chazelle}, Q.~{Jiu}, Q.~{Li}, and C.~{Wang}.
\newblock {Well-Posedness of the Limiting Equation of a Noisy Consensus Model
  in Opinion Dynamics}.
\newblock {\em ArXiv e-prints}, October 2015.

\bibitem{CPH:9676865}
Kai Jiang, Chu Wang, Yunqing Huang, and Pingwen Zhang.
\newblock Discovery of new metastable patterns in diblock copolymers.
\newblock {\em Communications in Computational Physics}, 14:443--460, 8 2013.

\end{thebibliography}
\bibliographystyle{unsrt}}

\vspace{3cm}
\section*{Appendix: Pseudo-Spectral Method for the PDE}

Pseudo-spectral method has been extensively adopted for simulating PDEs with periodic boundaries, especially in the fields of fluid dynamics and material science~\cite{fornberg1998practical,CPH:9676865}.
To solve Equation (\ref{PDE}) with periodic condition and $L=1$,
we expand $\rho(x,t)$ in Fourier space and cut off the expansion at $m$:

\begin{equation}
\rho(x,t)\approx\sum_{k=-m}^{m}\hat{\rho}_k(t)e^{i2\pi kx}.
\end{equation}
It follows that 
\begin{eqnarray*}
\int_{x-R}^{x+R}(y-x)\rho(y,t)\, dy 
&=& \int_{-R}^{R}y\rho(x+y,t)\, dy
\approx\sum_{k=-m}^{m}\hat{\rho}_k(t)e^{i2\pi kx}\int_{-R}^{R}ye^{i 2\pi ky}\, dy\\
&=&\sum_{-m\le k\le m, k\neq 0}\bigg\{\, \frac{R}{i2\pi k}\Bigl(e^{i2\pi kR}+e^{-i2\pi kR}\Bigr) \\
& & ~~~~~~~~~~~~~~~~~~~~
+\frac{1}{4\pi^2k^2}\Bigl(e^{i2\pi kR}-e^{-i2\pi kR}\Bigr)\bigg\}\hat{\rho}_k(t)e^{i2\pi kx}
\end{eqnarray*}

\noindent
Next we discretize the time as $\hat{\rho}_{k,r}=\hat{\rho}_{k}(rh)$
and generate the numerical scheme as follows: for all $-m\le k \le m$,

\begin{equation*}
\hat{\rho}_{k,r}\longrightarrow \bigg{[}\frac{R}{i2\pi k}(e^{i2\pi kR}+e^{-i2\pi kR})+\frac{1}{4\pi^2k^2}(e^{i2\pi kR}-e^{-i2\pi kR})\bigg{]}\hat{\rho}_{k,r}:=\hat{\varphi}_{k,r}
\end{equation*}
\begin{equation*}
\hat{\rho}_{k,r}\xrightarrow{FFT}\rho_{k,r}
\end{equation*}
\begin{equation*}
\hat{\varphi}_{k,r}\xrightarrow{FFT}\varphi_{k,r}
\end{equation*}
\begin{equation*}
\psi_{k,r}:=\varphi_{k,r}\rho_{k,r}
\end{equation*}
\begin{equation*}
\psi_{k,r}\xrightarrow{FFT}\hat{\psi}_{k,r}
\end{equation*}
\begin{equation*}
\hat{\rho}_{k,r+1}=\Big{(}-i2\pi k\hat{\psi}_{k,r}-2\pi^2\sigma^2 k^2\hat{\rho}_{k,r+1}\Big{)}h+\hat{\rho}_{k,r}
\end{equation*}

Notice that the constraint $\int \rho(y,t)\, dy=1$ is unconditionally satisfied in Equation~(\ref{PDE}); therefore we only need to set $\hat{\rho}_{0,r}=0$ during the iteration to prevent numerical error while no further treatment is required.

\end{document}